\documentclass{ifacconf}

\usepackage{graphicx}      % include this line if your document contains figures
\usepackage{natbib}        % required for bibliography

\usepackage{graphicx}      % include this line if your document contains figures
\usepackage{amssymb}
\usepackage{amsmath}
\usepackage{xcolor}
\usepackage{enumerate}
\usepackage{array}
\theoremstyle{definition}
\newtheorem{definition}{Definition }%[section]
\theoremstyle{plain}
\newtheorem{proposition}{Proposition }%[section]
\newtheorem{theorem}{Theorem }%[section]
\newcommand{\pfrac}[2]{\frac{\partial#1}{\partial #2}}
\newcommand{\dleq}[3]{#1 \leq #2 \leq #3}

\DeclareMathOperator{\mddiff}{\mathrm{ddiff}}

\newcommand{\mspan}{\mathrm{span}\,}
\newcommand{\mrk}{\mathrm{rk}\,}

\newcommand{\md}{\mathrm{d}}

\def\build#1_#2^#3{\mathrel{\mathop{\kern 0 pt#1}\limits_{#2}^{#3}}}

\def \vp {\varphi}
\def \a {\alpha}
\def \b {\beta}

\def \S {\Sigma}

\def \R {\mathbb{R}}

\def \t {\tilde}

\newcommand{\ol}[1]{\overline{#1}}
\begin{document}
\begin{frontmatter}

\title{Normal forms for $x$-flat   two-input  control-affine systems in dimension five}%\thanksref{footnoteinfo}} 

% \thanks[footnoteinfo]{Sponsor and financial support acknowledgment
% goes here. Paper titles should be written in uppercase and lowercase
% letters, not all uppercase.}

\thanks[footnoteinfo]{\begin{small}The second author has been supported by the Austrian Science Fund (FWF) under grant number P 32151.\end{small}}

\author[First,Second]{F. Nicolau}
\author[Third]{C. Gst{\"o}ttner}
\author[Fourth]{W. Respondek}

\address[First]{Quartz EA 7393, ENSEA, 95014 Cergy-Pontoise Cedex, France.}
\address[Second]{
Universit\'e Paris-Saclay, CNRS, CentraleSupélec, L2S, 91190, Gif-sur-Yvette, France, florentina.nicolau@ensea.fr.}
\address[Third]{Institute of Automatic Control and Control Systems Technology, Johannes Kepler University, Linz, Austria, conrad.gstoettner@jku.at.}
\address[Fourth]{Normandie Universit\'e, INSA Rouen, LMI, 76801 Saint-Etienne-du-Rouvray, France, witold.respondek@insa-rouen.fr.}

% \begin{abstract}  In this paper, we give normal forms for flat two-input control-{affine} systems in dimension five that admit a flat output depending on the state only (we call systems with that property $x$-flat systems). We discuss relations of $x$-flatness in dimension five  with static and dynamic feedback linearization and show that a system is $x$-flat if and only if it becomes linearizable via at most three prolongations of a suitably chosen control.  
% Therefore  $x$-flat systems in dimension five can be, in general, brought into normal forms generalizing the Brunovsk\'y canonical form. If  a system becomes linear via at most two-fold prolongation, the normal forms are structurally similar to the Brunovsk\'y form: they have a special triangular structure consisting  of a linear chain and a nonlinear one with at most two nonlinearities. If  a system becomes linear via a three-fold prolongation, we {obtain not only triangular structures but also a nontriangular one, and} face new interesting phenomena. \vspace{-0.4cm}

% \end{abstract}

\begin{abstract}  In this paper, we give normal forms for flat two-input control-{affine} systems in dimension five that admit a flat output depending on the state only (we call systems with that property $x$-flat systems). We discuss relations of $x$-flatness in dimension five  with static and dynamic feedback linearization and show that {if a system is $x$-flat} it becomes linearizable via at most three prolongations of a suitably chosen control.  
Therefore  $x$-flat systems in dimension five can be, in general, brought into normal forms generalizing the Brunovsk\'y canonical form. If  a system becomes linear via at most two-fold prolongation, the normal forms are structurally similar to the Brunovsk\'y form: they have a special triangular structure consisting  of a linear chain and a nonlinear one with at most two nonlinearities. If  a system becomes linear via a three-fold prolongation, we {obtain not only triangular structures but also a nontriangular one, and} face new interesting phenomena. \vspace{-0.4cm}

\end{abstract}


\end{frontmatter}
%===============================================================================

\section{Introduction}\vspace{-0.35cm}

In this paper, we give normal forms for flat control-affine systems of the form%\vspace{-0.05cm}
\begin{equation} \label{eq Sigma}
\Sigma\ :\ \dot x = f(x)+ g_1(x) u_1 + g_2(x) u_2,%\vspace{-0.05cm}
\end{equation}
where $x$ is the state defined on a open subset $X$ of~$\mathbb{R}^5$ {(more generally, a five-dimensional manifold~$ X$) and  $u = (u_1, u_2)$ is the control taking values in $\R^2$,} and where~$f$,~$g_1$ and~$g_2$ are {$\mathcal{C}^\infty$-smooth.}
%smooth. The word smooth will always mean $\mathcal{C}^\infty$-smooth.
The notion of flatness, {see Section~\ref{sec: flatness}}, was introduced in control theory in the {1990s}, by Fliess, L\'evine, Martin and Rouchon (\cite{fliess61vine},  see also {\cite{isidori1986sufficient, jakubczyk1993invariants, aranda1995linear}}, \cite{pomet1995differential}) and has attracted a considerable interest %{\cite{fliess1999lie,pomet1997dynamic,van1998differential,da2001relative}} 
because of its important applications in the problem of motion planning and constructive controllability.  
% (see, e.g., {\cite{martin2003flat,levine2009analysis, tang2011differential2,kolar2017time}}). 
%
%% REMOVED 1st submission
% {The} system  $\Xi: \dot x=~F(x,u)$, where $x \in X\subset \mathbb{R}^n$ and $u \in U\subset \mathbb{R}^m$, is \textit{flat} if we can find {locally}~$m$ functions $\varphi_i(x,u,\dots,u^{(r)})$, for some $r \geq 0$, such~that\vspace{-0.05cm}
% $$
% x=\gamma(\varphi,\dots,\varphi^{(s-1)}) \mbox{ and }
% u=\delta(\varphi,\dots,\varphi^{(s)}),%\vspace{-0.1cm}
% $$
% for a certain integer $s$ and suitable smooth maps $\gamma$ and~$\delta$, where $\varphi=(\varphi_1,\dots,\varphi_m)$ {is called a \textit{flat output}}. If all functions $\varphi_i$ depend on the state only, i.e., $\varphi_i =\varphi_i(x)$, for all $\dleq 1 i m$, we say that the system is $x$-flat. 
% % {Therefore, for a flat system, the evolution in time of all state and control variables can be recovered from that of flat outputs without integration and all trajectories of the system can be completely parameterized.}
%%%%%
%
{Systems} linearizable via invertible static feedback are flat and
their normal forms are well known:  they are static feedback equivalent to the Brunovsk\'y canonical form. In general, flat systems are not static feedback linearizable but  can be seen as a generalization of linear systems. Namely they are linearizable via dynamic, invertible and endogenous feedback, see \cite{fliess61vine,pomet1995differential,pomet1997dynamic}.  Although flatness has received a lot of attention because of its important applications, the problem of giving a catalogue of all 
normal forms for flat systems remains widely open (with some notable exceptions, like normal forms for the classes of static feedback linearizable systems \cite{brunovsky1970classification}, control-affine systems with $n-1$ controls \cite{charlet1989dynamic}, driftless systems with two controls \cite{martin1994feedback,murray1994nilpotent,li2010flat}, dynamically linearizable systems via one- or two-fold prologation of a suitably chosen control \cite{nicolau2019normal, nicolau2020normal}, {systems static feedback equivalent to a triangular form {\cite{silveira2015flat}, \cite{li2016multi} and} \cite{gstottner2021structurally,gstottner2022structurally}},   control-affine systems with two inputs and four states \cite{pomet1997dynamic}). 

The goal of this paper is to provide all normal forms for {all} $x$-flat {(i.e., admitting flat outputs depending on the state only)} two-input control-affine systems in dimension five {and thus can be considered as a continuation of \cite{pomet1997dynamic}}. 
Solving {that} problem in the simplest case of two controls and five states is interesting for few reasons; first, it yields a complete analysis for  a well defined class of flat systems, and second, it shows what kind of difficulties one may face when trying to give normal forms or to characterize flatness in the general case. 
We will prove that, in general,  $x$-flat control systems in dimension five can be brought into a normal form generalizing that of Brunovsk\'y. Indeed, $x$-flat control systems in dimension five  are either static feedback linearizable  or become static feedback linearizable after  prolonging at most three times a suitably chosen control. 
{For most cases, we obtain normal forms that are structurally similar to the Brunovsk\'y form: they have a special triangular structure consisting  of a linear chain and a nonlinear one with at most three nonlinearities but, in addition to that, we also face new interesting phenomena (for instance, {the lack of a triangular structure in}  {a particular case}).}
Therefore, an interesting observation is that even in a relatively small dimension, the most general triangular forms fail to capture all $x$-flat systems, and that illustrative normal forms do exist for all possible cases.
{Our results are related to those of \cite{silveira2007triangular}, where a triangular form for control-affine systems with two controls and at most  five states is presented. As expected, the triangular form of \cite{silveira2007triangular} in dimension five coincides with one of the forms of the catalogue presented in this paper.} 
%
% The paper is organized as follows. In Section~\ref{sec: flatness}, we recall the {definitions} of flatness and of differential weight. In Section~\ref{sec: main results}, 
% we give our main results
% {and {then} illustrate them} {an example} in Section~\ref{sec: examples}. We present  a sketch of proofs of our main theorems in Section~\ref{sec: proofs}.
\vspace{-0.25cm}

\section{Flatness} \label{sec: flatness} \vspace{-0.3cm}
Let $U$ be an open subset of $ \mathbb{R}^{ m}$. Fix an integer $l \geq -1$ and denote $U^l = U \times \mathbb{R}^{ml}$ and $\bar u^l = (u,\dot  u, \dots, u^{(l)})$. For $l=-1$,  the set $ U^{-1} $ is empty and~$\bar u^{-1}$ in an empty sequence.

\begin{definition} \label{Def platitude}
The system $\Xi: \dot x=~F(x,u)$, $x\in X\subset \mathbb{R}^n$, $u\in U\subset \mathbb{R}^m$, is \textit{flat} at $(x_0, \bar u_0^l) \in X \times U^l$, for  $l\geq -1$, if there exist a neighborhood~$\mathcal O^{l}$ of $(x_0, \bar u_0^l)$
and~$m$ smooth functions $\varphi_i=\varphi_i(x, u,\dot  u, \dots, u^{(l)})$, $1 \leq i\leq m$, defined in~$\mathcal O^l$, having the following property: there exist an integer~$s$ and smooth functions $ \gamma_i$, $1 \leq i \leq n$, and $\delta_j$, $ 1 \leq j \leq m$, such that \vspace{-0.05cm}
$$
x_i = \gamma_i (\varphi, \dot \varphi, \dots, \varphi^{(s-1)}) \mbox{ and }u_j=  \delta_j (\varphi, \dot \varphi, \dots, \varphi^{(s)})
$$
for any $C^{l+s}$-control $u(t)$ and corresponding trajectory $x(t)$
that satisfy $(x(t), u(t), \dots ,$ $u^{(l)}(t)) \in~\mathcal O^l$, where $\varphi=(\varphi_1, \dots, \varphi_m)$ and is  called a \textit{flat output}. {If $\varphi_i =\varphi_i(x)$, for all $\dleq 1 i m$, we say that the system is $x$-flat.}
\end{definition}

The minimal number of derivatives of components of a flat output, needed to express~$x$ and~$u$, will be called the  differential weight of that flat output and is formalized as follows.
By definition, for  any flat output $\varphi$ of $\Xi$ there exist integers $s_1,\dots,s_m$ such that\vspace{-0.05cm}
\begin{equation}\label{eq:descript}
\left.
\begin{array}{lll}
x&= &\gamma(\varphi_1,\dot\varphi_1,\dots,\varphi_1^{({s_1-1})},\dots, \varphi_m,\dot\varphi_m, \dots, \varphi_m^{({s_m-1})}) \\
u&= & \delta(\varphi_1,\dot\varphi_1,\dots,\varphi_1^{(s_1)},\dots, \varphi_m,\dot\varphi_m, \dots, \varphi_m^{(s_m)}).
\end{array}
\right.\vspace{-0.05cm}
\end{equation}
Moreover, we can choose $(s_1,\dots,s_m)$, {$\gamma$ and $\delta$} such that (see \cite{respondek2003symmetries}) if for any other~$m$-tuple $(\tilde s_1,\dots,\tilde s_m)$ and functions {$\t \gamma$ and $\t \delta$}, we have\vspace{-0.05cm}
$$
\left.
\begin{array}{lll}
x&= &\tilde \gamma(\varphi_1,\dot\varphi_1,\dots,\varphi_1^{({\tilde s_1-1})},\dots, \varphi_m,\dot\varphi_m, \dots, \varphi_m^{({\tilde s_m-1})}) \\
u&= &\tilde \delta(\varphi_1,\dot\varphi_1,\dots,\varphi_1^{(\tilde s_1)},\dots, \varphi_m,\dot\varphi_m, \dots, \varphi_m^{(\tilde s_m)}),
\end{array}
\right.
$$
then $s_i \leq \tilde s_i$, for $\dleq 1 i m$. We will call  $\sum_{i=1}^{m}(s_i+1) =m+\sum_{i=1}^{m}s_i$ the differential weight of $\varphi$. A flat output of $\Xi$ is called \textit{minimal} if its differential weight is the lowest among all flat outputs of $\Xi$. We define the \textit{differential weight} of a flat system to be equal to the differential weight of a minimal flat output.
The differential weight {is} $n+m+{p}$, where ${p}\geq 0$ {and} can be interpreted as the minimal dimension of a precompensator that dynamically linearizes the system, {and we will call it the {\it differential difference}, see \cite{gstottner2021finite, gstottner2021necessary}, shortly, the ddiff of a flat system}.

Consider two  control-affine {systems} $\Sigma$ and $\t \Sigma$ given, resp., by $\dot  x = f(x)+\sum_{i=1}^mu_ig_i(x)$, $x\in X$, $u\in {\mathbb{R}^m}$ and $\dot {\t x} =\t  f(\t x)+\sum_{i=1}^m\t u_i\t g_i(\t x)$, $\t x\in \t X$, $\t u\in {\mathbb{R}^m}$, where $X$ and $\t X$ are open subsets of $\mathbb{R}^n$. We say that  $\Sigma$ and $\t \Sigma$ are {\it locally static feedback equivalent} if there exist a (local) diffeomorphism $\t x = \phi(x)$ and an invertible static feedback transformation  {of the form} $u=\alpha(x)+\beta(x)\t u $ which transform the first system into the second, i.e., $\t f(\phi(x)) = \pfrac{\phi(x)}{x}(f(x)+g(x) \a(x))$ and  $\t g(\phi(x)) = \pfrac{\phi(x)}{x} g(x)\b(x)$, where $g = (g_1, \ldots, g_m)$ and $\t g = (\t g_1, \ldots, \t g_m)$. 
The control-affine system $\Sigma$
is {\it static feedback linearizable} if it is static feedback  equivalent to a linear controllable system {$\Lambda : \dot {\t x}=A\t x+B\t u$}. The problem of static feedback linearization was solved {by} 
%\cite{brockett1979feedback} (for a smaller class of transformations) and then by 
\cite{jakubczyk1980on} and, independently, by \cite{hunt1981linear}, who gave the following geometric necessary and sufficient conditions.
Define
the distributions $\mathcal{D}^{ j+1} = \mathcal{D}^{ j} + [f,\mathcal{D}^{ j}]$, where $\mathcal{D}^0 = \mbox{span} \{g_1,\ldots, g_m\}$ and $[f,\mathcal{D}^{ j}] =\{[f,\xi]: \xi \in \mathcal{D}^{ j}\}$.  The system $\Sigma$ is  locally static feedback linearizable if and only if for any $ j \geq 0$, the distributions~$\mathcal{D}^{ j}$ are of constant rank, \textit{involutive} and $\mathcal{D}^{n-1} = TX$.
% {Therefore, the following {nested} sequence of  involutive distributions summarizes the geometry of static feedback linearizable systems:
% $
% \mathcal{D}^0 \subset \mathcal{D}^1\subset \ldots \subset\mathcal{D}^{n-1} = TX.
% $}
A feedback linearizable system is static feedback equivalent to the Brunovsk\'y canonical form \cite{brunovsky1970classification}
\vspace{-0.1cm}
$$
(Br):\left\{\quad
\begin{array}{lcl}
\dot z_{i}^{j} & = & z_{i}^{j+1}\\
\dot z_{i}^{\rho_i} & = & v_i,
\end{array}\right.\vspace{-0.1cm}
$$
where $1\leq i\leq m$, $1\leq j\leq \rho_i-1$, and $\sum_{i=1}^m \rho_i =n$,
and is clearly flat with  $\varphi = (\vp_1, \ldots, \vp_m) =$ $(z_{1}^{1}, \ldots, z_{m}^{1})$ being a minimal flat output (of differential weight $n+m$, {i.e., ddiff{$(Br)$} = 0}).
In fact, an equivalent way of describing static feedback linearizable systems is that they are flat systems of differential weight {$n+m$}, {equivalently, of ddiff 0}, see Theorem~2.2 in \cite{nicolau2017flatness}.

The goal of this paper is to provide normal forms for two-input $x$-flat control systems in dimension five. We will prove that, in general,  $x$-flat control systems in dimension five can be brought into a normal form generalizing that of Brunovsk\'y. Namely, we will show that they are either static feedback linearizable  or become static feedback linearizable after  prolonging at most three times a suitably chosen control, which is the simplest dynamic feedback. 
%
%In the {former} case, they are static feedback equivalent to the  Brunovsk\'y canonical form. 
% In the later case, at most three nonlinearities can appear in the normal forms. If the system becomes linear via one- or two-fold prolongation, the normal forms are structurally similar to the Brunovsk\'y canonical form. If the system becomes linear via a three-fold prolongation, we 
% {obtain not only triangular structures but also a nontriangular one and new  phenomena appear (that will be discussed in details in Section~\ref{sec: discussion}).} 
%
{For} systems that are not static feedback linearizable, 
(not flat of differential weight $n+m$),
there exists {the smallest} integer $\dleq 0 k n-1$ such that the linearizability conditions (either involutivity or constant rank) are not satisfied for $\mathcal{D}^k$ (flat systems are always accessible so $\mathcal D^{n-1}=TX$ holds). 
We work under constant rank assumptions (see 
Assumption 
{(A)} at the begining of Section~\ref{sec: main results}), {so only the case of~$\mathcal{D}^k$} noninvolutive can occur. 
We will see when discussing the normal forms that the {integer $k$ indicating} the first noninvolutive disctribution plays an important role: namely, for flatness singularities in the control space. \vspace{-0.35cm}

\section{Main results: normal forms}\vspace{-0.3cm}
\label{sec: main results}

We study $x$-flat  two-input control-affine systems in dimension five {of the form}:\vspace{-0.1cm}
\begin{equation} \label{Sigma}
\Sigma\ :\ \dot x = f(x)+u_1g_1(x)+u_2g_2(x),\vspace{-0.05cm}
\end{equation}
where $x\in X$, an open subset of $\mathbb{R}^5$, and~$u\in \mathbb{R}^2$.  
{Throughout we make the following} {assumption:}\vspace{-0.05cm}
%{The proofs of all results of this paper will be provided in a future publication. We make the following assumption:}

{\bf Assumption {(A)}} {\it 
 All ranks involved are constant in a neighborhood of a given 
${(x_0,  \bar u_0^{l})} \in X \times \mathbb{R}^{2(l+1)}$, for a certain  $\dleq {-1}l 1$. 
All results are thus valid on an open and dense subset of $ X \times  \mathbb{R}^{2(l+1)}$
and hold locally, around any given point of that set.
% From now on, unless stated otherwise, we assume that all ranks involved are constant in a neighborhood of a given 
% ${(x_0,  \bar u_0^{l})} \in X \times \mathbb{R}^{2(l+1)}$, for a certain  $\dleq {-1}l 1$. 
% All results presented here are thus valid on an open and dense subset of $ X \times  \mathbb{R}^{2(l+1)}$
% and hold locally, around any given point of that set.
}

The following result {relates} $x$-flatness in dimension five with differential weight and {with} dynamic feedback linearization via prolongations of a suitably chosen control. {Its proof, and those of all results of this paper will be provided in a future publication.}

\begin{proposition}\label{prop: x flatness dw linearization}
{\it The following {statements} are equivalent:\vspace{-0.15cm}
   \begin{enumerate}[\normalfont(i)]
    \item $\Sigma $ is $x$-flat at $(x_0, u_0,\dot  u_0, \dots, u_0^{(l)})$.
    \item $\Sigma $ is $x$-flat, at either~$x_0$ or $(x_0,u_0)$ or $(x_0,u_0, \dot u_0)$, of differential weight at most  $n+m+ 3 = 10$, {that is, $\mathrm{ddiff} (\S)\leq 3$}.
        \item There exists, around~$x_0$, an invertible static feedback transformation $u = \alpha(x) + \beta(x)\tilde u,$ bringing~$\Sigma$  into the form $\tilde \Sigma : \dot x =\tilde  f(x)+\tilde  u_1 \tilde  g_1(x) +\tilde   u_2\tilde  g_2(x)$, such that the prolonged system
\vspace{-0.15cm}
$$\tilde {\Sigma}^{(p, 0)} :
\left \{
\begin{array}{lcl}
\dot x &= &\tilde f(x)+ y_1 \tilde g_1(x)+v_2\tilde  g_2(x), \\
\dot y_i&= &y_{i+1}, \quad {\dleq 1 i p-1,}\\
\dot y_p&= &v_1,
\end{array}
\right.\vspace{-0.15cm}
$$
%where $\dleq 1 i p-1$, 
for a certain $\dleq 0 p 3$, is locally static feedback linearizable {{around} $(x_0,y_{0})$}{, with a linearizing output depending on $x$ only\footnote{This is always the case if $p\leq 2$ and the condition is needed for $p=3$ only.}, and} with $y_1 = \tilde u_1$, $v_2 = \tilde u_2$,  $\tilde f = f +g  \alpha $ and  $\tilde g = g\beta$, where $g = (g_1, g_2)$ and $\tilde g = (\tilde g_1,\tilde g_2)$. 
\end{enumerate}\vspace{-0.15cm}
{Moreover, the $\mddiff(\S)$ of item (ii) coincides with the smallest integer $p$ of item (iii).}
}
\end{proposition}

%% Without ``enumerate''
%
% \begin{proposition}\label{prop: x flatness dw linearization}
% {\it The following are equivalent:
% %    \begin{enumerate}[(i)]
% %     \item 
% 
% (i) \quad $\Sigma $ is $x$-flat at $(x_0, u_0,\dot  u_0, \dots, u_0^{(l)})$.
% 
% (ii) \quad $\Sigma $ is $x$-flat, at either~$x_0$ or $(x_0,u_0)$ or $(x_0,u_0, \dot u_0)$, of differential weight at most  $n+m+ 3 = 10$, {that is, $\mathrm{ddiff} (\S)\leq 3$}.
%     
%     {Q: $ddiff < 3$ or $\mathrm{ddiff} < 3$, i.e., in italic or not?}
% 
% (iii) \quad There exists, around~$x_0$, an invertible static feedback transformation $u = \alpha(x) + \beta(x)\tilde u,$ bringing~$\Sigma$  into the form $\tilde \Sigma : \dot x =\tilde  f(x)+\tilde  u_1 \tilde  g_1(x) +\tilde   u_2\tilde  g_2(x)$, such that the prolonged system
% \vspace{-0.1cm}
% $$\tilde {\Sigma}^{(p, 0)} :
% \left \{
% \begin{array}{lcl}
% \dot x &= &\tilde f(x)+ y_1 \tilde g_1(x)+v_2\tilde  g_2(x), \\
% \dot y_i&= &y_{i+1}, \quad \\
% \dot y_p&= &v_1,
% \end{array}
% \right.\vspace{-0.1cm}
% $$
% where $\dleq 1 i p-1$, for a certain $\dleq 0 p 3$, 
% is locally static feedback linearizable {{around} $(x_0,y_{0})$}, with $y_1 = \tilde u_1$, $v_2 = \tilde u_2$,  $\tilde f = f +g  \alpha $ and  $\tilde g = g\beta$, where $g = (g_1, g_2)$ and $\tilde g = (\tilde g_1,\tilde g_2)$.
% %\end{enumerate}
% 
% {Moreover, the $\mddiff(\S)$ of item (ii) coincides with the smallest integer $p$ of item (iii).}
% }
% \end{proposition}

A system $\Sigma$ satisfying (iii) with $p=0$ is simply static feedback linearizable. If  $\Sigma$ satisfies  (iii) with $\dleq 1 p~3$, then it  will be called \emph{dynamically linearizable via {an} invertible $p$-fold prolongation}. $\tilde \Sigma^{(p, 0)}$ is, indeed, obtained by applying an invertible static feedback $u = \alpha + \beta\tilde u$ and then prolonging  the first control $\tilde u_1$ $p$-times as $v_1 =  {\tilde u}_1^{(p)}$ and not prolonging $\tilde u_2$ (which explains the notation $\tilde \Sigma^{(p, 0)}$). 
It follows that any two-input $x$-flat system in dimension five is of differential weight at most $n+m+3 =10 $ and becomes static feedback linearizable after an at most three-fold prolongation of a suitably chosen control. 

The main results are given by Theorems~\ref{thm: normal forms dw at most 9} and \ref{thm: normal forms dw 10} that present several normal forms  for the class of flat two-input control-affine  systems in dimension 5. 
Given the system~$\S$ that is $x$-flat around $x_0$ or around $(x_0,u_0)$  or around $(x_0,u_0, \dot u_0)$, {all} normal forms are obtained under the local static feedback transformation\vspace{-0.1cm}
$$
z= \phi(x), \qquad u = \a(x) + \b(x) v,\vspace{-0.1cm}
$$
and are $x$-flat around  $z_0 =  \phi(x_0)$ or around $(z_0, v_0)$  or around $(z_0, v_0, \dot v_0)$, where $ u_0 = \a(x_0) + \b(x_0)v_0$ and 
{$\dot u_0 = L_F(\a + \b v_0)(x_0, u_0) +  \b(x_0)\dot v_0 $, with $F = f + gu$.}
%$\dot u_0 =\dot  \a(x_0, u_0) + \dot \b(x_0, u_0)v_0 +  \b(x_0)\dot v_0 $. 
Below,~$v_{10}$ stands for the first component of $v_0$.  
For $\dleq 1 i  5$, $\ol z_i$ denotes  $\ol z_i = (z_1, z_2, \ldots, z_i)$.

\begin{theorem}\label{thm: normal forms dw at most 9}
{\it The following  {statements}  are equivalent: \vspace{-0.1cm}
 \begin{enumerate}[\normalfont (i)]
  \item  $\Sigma $ is $x$-flat at either~$x_0$ or $(x_0,u_0)$, of differential weight at most  $n+m+ 2 = 9$, {i.e., $\mddiff(\S) \leq 2$}.
  \item  $\Sigma $ is locally, around~$x_0$, static feedback equivalent in a neighborhood of $z_0\in \R^n$
either to\vspace{-0.05cm} 
$$
NF: \left\{\begin{array}{rcl rcl}
  \dot z_1& = &v_1 & \dot z_2& = & z_3, \\
  & & & \dot z_3& = &b_1(\ol z_4) + b_2(\ol z_4) v_1, \\%& \pfrac{(b_1 +b_2 v_{10})}{z_4}(z_0)\neq 0,   \\
  & & & \dot z_4& = &a_1(z) + a_2(z)v_1, \\ %& \pfrac{(a_1 +a_2 v_{10})}{z_5}(z_0)\neq 0,   \\
  & & & \dot z_5& = &v_2,
 \end{array}\right.\vspace{-0.05cm} 
 $$
 or to \vspace{-0.05cm} 
  $$
NF':  \left\{ \begin{array}{r@{}c@{}l rcl}
  \dot z_1& = &z_2 & \dot z_3& = &b_1(\ol z_4) + b_2(\ol z_4) v_1,\quad  b_2(z_0) = 0, \\
  \dot z_2& = &v_1& \dot z_4& = &a_1(z) + a_2(z)v_1,\\ %& \pfrac{(a_1 + a_2 v_{10})}{z_5}(z_0)\neq 0, \\
  & & & \dot z_5& = &v_2,  
 \end{array}\right.\vspace{-0.05cm} 
 $$
 %
 %{Q: use $d_3$ and $d_4$ instead of $a$ and $b$ (see relation (5))?}
{with  $ \pfrac{{b}}{z_4}(z_0, v_{10})\neq 0$ and  $ \pfrac{{a}}{z_5}(z_0, v_{10})\neq 0$, where $ {b} = b_1 +b_2 v_{1}$ and  $ {a}=   a_1 +a_2 v_{1} $,}
  %where $ \pfrac{(b_1 +b_2 v_{10})}{z_4}(z_0)\neq 0$ and  $ \pfrac{(a_1 +a_2 v_{10})}{z_5}(z_0)\neq 0$
  for both $NF$ and $NF'$, 
or to\vspace{-0.05cm} 
  $$
NF_{ 7}:  \left\{ \begin{array}{r@{}c@{}l rcl}
  \dot z_1& = &z_2 & \dot z_4& = &a_1(z) +  (z_5-z_{50})v_1,\\
  \dot z_2& = &z_3& \dot z_5& = &v_2,\\
  \dot z_3& = &v_1,  
 \end{array}\right.\vspace{-0.05cm} 
 $$
where  $ \pfrac{a_1}{z_5}(z_0) + v_{10}\neq 0$. 
  \item  $\Sigma $ is locally, around~$x_0$, static feedback equivalent to one of the normal forms $NF^*_1$- $NF^*_{ 7}$, with $NF^*_1$-$NF^*_{ 6}$  detailed  in Table~\ref{table forms NF NF'}, and where $NF^*$ stands either for $NF$ or for $NF'$,
  {and $NF_7^*$ denotes $NF_7$.} %. {Q:  $NF_7$ OR $NF_7^*$?}
  
  \item $\Sigma $ is either static feedback linearizable or dynamically linearizable via an invertible one-fold prolongation or  dynamically linearizable via an invertible two-fold prolongation. 
 \end{enumerate}
}
\end{theorem}

The normal forms $NF$ and $NF'$ present four nonlinearities, but we can always normalize at least two of them  leading, respectively, to one of the normal forms $NF_1$-$NF_{5}$ {associated to $NF$ or to $NF'_1$-$NF'_{6}$ associated to $NF'$, see} Table~\ref{table forms NF NF'}. 
%(where d.w. stands for differential weight). %presents all possible cases: forms $NF_1$-$NF_{5}$ are normalizations associated to $NF$, while  $NF'_1$-$NF'_{6}$  are associated to $NF'$.  
% When  we say that $NF'_{5}$ is missing, we mean that there is no form corresponding to $NF'$ with the same normalizations as those leading to $NF_5$ for $NF$ (such a normalization for  $NF'$  is actually static feedback equivalent to $NF'_1$). A similar remark can be done for $NF_6$. 
Notice that there is no form {$NF_6$} corresponding to $NF$ with the same normalizations as those leading to $NF'_6$ for $NF'$ (such a normalization for  $NF$  is actually static feedback equivalent to $NF'_1$).
Like in item (iii) of Theorem~\ref{thm: normal forms dw at most 9}, we will use the notation $NF^*_1$-$NF^*_{ 7}$, where $NF^*$ stands either for $NF$ or for $NF'$, form $NF^*_{ 7}$  refers to $NF_7$ of  Theorem~\ref{thm: normal forms dw at most 9}~(ii), while $NF^*_1$- $NF^*_{ 6}$ denote the forms of  Table~\ref{table forms NF NF'}. 
Conditions in the column ``Input regularity conditions'' of Table~\ref{table forms NF NF'} describe pairs $(z_0, v_{10})$ around which the system is flat{; notice that the ``Input regularity conditions'' mix up state and controls, implying that around any $z_0$ there are controls that destroy flatness: for instance, for $NF^*_2$, we have to avoid $\pfrac{a_1}{z_5}(z) + v_1 =0$ (which has a solution in any neighborhood of %a ``good''
any fixed $z_0$).}
For some of the forms, the nonlinear functions should additionally satisfy some conditions such {as}:
\\
{\bf (C1)} {\it~$\S$ is not static feedback linearizable (i.e., it does not satisfy the conditions of \cite{jakubczyk1980on}).}
\\
or 
\\
{\bf  (C2)} {\it~$\S$ is not dynamically linearizable via an invertible one-fold prolongation (i.e., it does not satisfy the conditions of \cite{nicolau2017flatness}).
}

% \begin{remark}
% \label{rk: struct cond}
Conditions (C1) and (C2) are structural conditions that imply the noninvolutivity of some distributions (see Section~\ref{sec: flatness}). For the nonlinear functions of the normal forms, {noninvolutivity} translates into derivatives with respect to certain variables {being} not identically zero, so it is expressed by non-equalities and, in this sens, {(C1) and (C2)}  are regularity conditions (thus we call them structural regularity conditions in the tables presenting our normal forms).
%\end{remark}
%
%
For instance, if (C2) is not {assumed} for $NF_5$, then the case $\pfrac{b_2}{z_4}\equiv 0$ could occur, which would imply local static feedback equivalence to either $NF'_1$ or $NF'_6$.
%%%%%%%REMOVED 1st submission
% Without (C2), $NF_5$ could also be locally static feedback equivalent to $NF_1$ ({if} $b_2\equiv 0$). $NF_5$ can never be equivalent to $NF^*_2$ or $NF^*_3$ because of the distribution $\mathcal D^0$ being involutive.
% Moreover, the function~$b_2$ of $NF_5$ is assumed to be such that $\pfrac{b_2}{z_4}(z_0)=0$ (if $\pfrac{b_2}{z_4}(z_0)\neq 0$, form   $NF_5$ is locally static feedback equivalent to $NF_4$). 
%%%%%%%%%%%%%%%

{Notice} the two different normalizations  $b_2 = z_4$ for $NF_3$ and  $NF_4$,  and   $b_2 = z_4 -z_{40}$ for $NF'_3$ and $NF'_4$. This is a consequence of the fact that the function $b_2$ associated to $NF'$ has to satisfy $b(z_0) = 0$; otherwise, $NF$ and $NF'$ are actually static feedback equivalent (indeed by applying the invertible feedback transformation  $\t v_1 = b_1 + b_2 v_1$, $NF'$ would be transformed into $NF$). {A similar remark} also applies for $NF_7$. 

\begin{table*}%[hb]
\begin{center}
\caption{Simplification {of}  
%the normal forms 
$NF$ and $NF'$ {(with d.w. standing for differential weight)}.}\label{table forms NF NF'}
\begin{tabular}{>{\raggedright}m{0.5cm}>{\raggedright}m{1.7cm}>{\raggedright}m{3.1cm}>{\raggedright}m{5cm}>{\raggedright}m{2.5cm}m{2.8cm}}
ddiff & d.w. & Forms & Nonlinearities  & {Input regularity conditions} &{Structural regularity conditions}\\
\hline
0& $n+m=7$ & $NF_1$, $NF'_1$ &  $b_1 = z_4$, $b_2\equiv 0$ &\hspace{1cm} -& \hspace{1.3cm} - \\ 
&& &   $a_1 = z_5$, $a_2\equiv 0$ &\hspace{1cm} -& \hspace{1.3cm} - \\ 
\hline 
1& $n+m+1=8$  & $NF_2$, $NF'_2$ &   $b_1 = z_4$, $b_2\equiv 0$ &\hspace{1cm} -& \hspace{1.3cm} - \\ 
&& &  $a_1(z)$ any, $a_2 = z_5$ &  $\pfrac{a_1}{z_5}(z_0)+ v_{10}\neq 0$  &  \hspace{1.3cm} -\\  
\hline
2& $n+m+2=9$  & $NF_3$,  $NF'_3$ & $  b_1 (\ol z_4)$ any, $b_2 = \left\{\begin{array}{l}
                                                                    z_4, \mbox{ for } NF_3\\
                                                                    z_4-z_{40}, \mbox{ for }  NF'_3
                                                                   \end{array}\right.$
     &  $\pfrac{b_1}{z_4}(z_0)+ v_{10}\neq 0$&\hspace{1.3cm} -\\  
& &&   $a_1(z)$ any, $a_2 = z_5$ & $\pfrac{a_1}{z_5}(z_0)+ v_{10}\neq 0$& \hspace{1.3cm} \\ 
\hline
2&$n+m+2=9$  & $NF_4$, $NF'_4$ & $  b_1  (\ol z_4)$ any, $b_2 = \left\{\begin{array}{l}
                                                                    z_4, \mbox{ for } NF_4\\
                                                                    z_4-z_{40}, \mbox{ for }  NF'_4
                                                                   \end{array}\right.$ 
                                                                   
     &  $\pfrac{b_1}{z_4}(z_0)+ v_{10}\neq 0$&\hspace{1.3cm} - %{needed for NF'? for NF no} 
     \\  
& &&   $a_1 = z_5$, $ a_2 (\ol z_4)$ any& \hspace{1cm}   - & \\ 
\hline
2&$n+m+2=9$  & $NF_5$, $NF'_5$ & $b_1= z_4$, $b_2 (\ol  z_4)$ any&\hspace{1cm}   -&  {(C2)} and $ \pfrac{b_2}{z_4}(z_0) = 0$ \\  
&& &   $a_1 = z_5$, $a_2 (\ol z_4)$ any&\hspace{1cm}  -& \\ 
\hline
1&$n+m+1=8$ &  $NF'_6$ (no corresp. $NF_6$) &   $b_1 = z_4$, $b_2\equiv 0$ &\hspace{1cm} -& \hspace{1.3cm} - \\ 
&& &   $a_1 = z_5$,  $a_2(\ol z_4)$ any& \hspace{1cm} - & \hspace{1.15cm}{(C1)} \\
\hline
\end{tabular}
\end{center}
\end{table*}

Theorem~\ref{thm: normal forms dw 10} {below} treats the case of differential weight $n+m+ 3 = 10$, {i.e., $\mddiff(\S)  = 3$},  and the nonlinear functions appearing in item (ii) {below} have to verify some additional {structural} regularity conditions: %such that: 
\\
{\bf (C3)} {\it $\Sigma $ is not linearizable via an invertible two-fold prolongation (i.e., it does not satisfy the conditions of \cite{nicolau2016flatness,gstottner2021necessary}).}

\begin{theorem}\label{thm: normal forms dw 10}
 {\it The following  {statements}  are equivalent: 
 \begin{enumerate}[\normalfont(i)]
  \item  $\Sigma $ is $x$-flat, at either  $(x_0,u_0)$ or  $(x_0,u_0, \dot u_0)$, of differential weight  $n+m+ 3 = 10$, {i.e., $\mddiff(\S)  = 3$.}
  \item  $\Sigma $ is locally, around~$x_0$, static feedback equivalent in a neighborhood of $z_0\in \R^n$
either to\vspace{-0.05cm}
$$
{NF''}: \left\{\begin{array}{rcl rcl}
  \dot z_1& = &v_1 & \dot z_2& = &c_1(\ol z_3) + c_2(\ol z_3)v_1 , \\%& \pfrac{(c_1 +c_2 v_{10})}{z_3}(z_0)\neq 0, \\
  & & & \dot z_3& = &b_1(\ol z_4) + b_2(\ol z_4) v_1,\\%& \pfrac{(b_1 +b_2 v_{10})}{z_4}(z_0)\neq 0,   \\
  & & & \dot z_4& = &a_1(z) + a_2(z)v_1,\\%&\pfrac{(a_1 +a_2 v_{10})}{z_5}(z_0)\neq 0,   \\
  & & & \dot z_5& = &v_2,
 \end{array}\right.\vspace{-0.05cm}
$$
 with $ \pfrac{(c_1 +c_2 v_{10})}{z_3}(z_0)\neq 0$, $\pfrac{(b_1 +b_2 v_{10})}{z_4}(z_0)\neq 0$, $\pfrac{(a_1 +a_2 v_{10})}{z_5}(z_0)\neq 0$, and 
 the functions $a_i, b_i, c_i$, $\dleq 1 i 2$, verifying {the} additional regularity  conditions {that guarantee condition (C3) to hold},
or to \vspace{-0.05cm}
$$
 NF_{ {13}}:\left\{ \begin{array}{rcl rcl}
  \dot z_1& = &v_1 & \dot z_2& = &z_3 + z_4v_1 \\
  & & & \dot {z}_3& = &a(\ol z_4) + (-z_5 + b(\ol z_4)) v_1  \\
  & & & \dot {z}_4& = &z_5 + c(\ol z_4)  v_1 \\
  & & & \dot z_5& = &v_2, 
  \end{array}\right.\vspace{-0.05cm}
$$
where \vspace{-0.05cm}
\begin{equation} \label{eq: reg cond}
\begin{array}{c}
    \left(\pfrac{a}{z_4} - (\pfrac{a}{z_3} - \pfrac{b}{z_4})v_{10} -  (\pfrac{b}{z_3} - \pfrac{c}{z_4})v^2_{10} \right.\\
   \left. -\pfrac{c}{z_3}v^3_{10}+ \dot v_{10}  \right)(z_0)\neq 0.
    \end{array}
\end{equation}

  \item  $\Sigma $ is locally, around~$x_0$, static feedback equivalent to one of the normal forms $NF''_8$-$NF''_{12}$ or $NF_{ {13}}$, with $NF''_8$-$NF''_{12}$ detailed  in Table~\ref{table form NF''}. 
  
  %\item $\Sigma $ is  dynamically linearizable via an invertible three-fold prolongation. 
 \end{enumerate}
{Moreover, if any of the equivalent items $(i)\Leftrightarrow(ii)\Leftrightarrow(iii)$ holds, then the system is linearizable by an invertible three-fold prolongation.}
}
\end{theorem}
{The} normal form $NF''$ presents six nonlinearities but we can always normalize at least one nonlinearity per nonlinear equation, leading to one of the normal forms $NF''_8$-$NF''_{12}$  presented in Table~\ref{table form NF''} below. The non normalized nonlinearities are such that {condition (C3)} is satisfied for all forms (all of them being of differential weight $n+m+3 = 10$); {some of them may depend explicitly on some $z$'s but their derivatives with respect to those $z$-variables have to vanish at~$z_0$} (see, for instance, the function~$b_2$ of $NF''_{10}$).\vspace{-0.2cm}

\begin{table*}%[h]
\begin{center}
\caption{Simplification {of}
%the normal form 
{$NF''$} ($\mathrm{d.w.}(\S)  = n+m+3=10$ and $\mddiff(\S) = 3$).}
%,  {with d.w. standing for differential weight}).}
\label{table form NF''}
 \begin{tabular}{>{\raggedright}m{1cm}>{\raggedright}m{3.4cm}>{\raggedright}m{3.5cm}m{4.6cm}}
 Form & Nonlinearities  &{Input regularity conditions} &{Structural regularity conditions}\\
  \hline
   $NF''_8$ &  $c_1(\ol  z_3)$ any, $c_2 = z_3$ &   $\pfrac{c_1}{z_3}(z_0)+ v_{10}\neq 0$& \hspace{2cm} - \\ 
    &  $b_1(\ol  z_4)$ any, $b_2 = z_4$ &  $\pfrac{b_1}{z_4}(z_0)+ v_{10}\neq 0$& \\ 
  &   $a_1(z) $ any, $a_2= z_5$ & $\pfrac{a_1}{z_5}(z_0)+ v_{10}\neq 0$& \\ 
\hline
  $NF''_{9}$ &   $c_1(\ol  z_3)$ any, $c_2 = z_3$ &   $\pfrac{c_1}{z_3}(z_0)+ v_{10}\neq 0$&\hspace{2cm}{(C3)} \\ 
    &  $b_1(\ol  z_4)$ any, $b_2 = z_4$ &  $\pfrac{b_1}{z_4}(z_0)+ v_{10}\neq 0$ &\\ 
      &   $a_1 = z_5 $, $a_2(\ol  z_4)$ any &\hspace{1.4cm}  -& \\ 
\hline
 $NF''_{10}$ & $c_1(\ol  z_3)$ any, $c_2 = z_3$ &   $\pfrac{c_1}{z_3}(z_0)+ v_{10}\neq 0$&\hspace{1cm}{(C3)} and \\ 
    &  $b_1= z_4$, $b_2(\ol  z_4)$ any &\hspace{1.4cm} -  & \hspace{2.3cm} $\pfrac{b_2}{z_4}(z_0) = 0 $\\ 
    %$1+\pfrac{b_2}{z_4}(z_0)v_{10}\neq 0$ ,\pfrac{b_2}{z_4}\not\equiv0
      &   $a_1 = z_5 $, $a_2(\ol  z_4)$ any & \hspace{1.4cm} -& \\ 
\hline
 $NF''_{11}$ & $c_1 = z_3$, $c_2(\ol  z_3)$ any & \hspace{1.4cm}  -&\hspace{1cm} {(C3)}  and  $\pfrac{c_2}{z_3}(z_0) = 0$\\ 
 % $1+ \pfrac{c_2}{z_3}(z_0)v_{10}\neq 0$ , \pfrac{c_2}{z_3}\not\equiv0$
    &  $b_1(\ol  z_4)$ any, $b_2= z_4$ & $\pfrac{b_1}{z_4}(z_0)+ v_{10}\neq 0$ & \\ 
      &   $a_1 = z_5 $, $a_2(\ol  z_4)$ any &\hspace{1.4cm}  - &\\ 
\hline
 $NF''_{12}$ & $c_1 = z_3$, $c_2(\ol  z_3)$ any & \hspace{1.4cm}  - &\hspace{1cm} {(C3)}  and  $\pfrac{c_2}{z_3}(z_0) =  0$ \\  
 %,\pfrac{c_2}{z_3} \not\equiv0 
    &  $b_1= z_4$, $b_2(\ol  z_4)$ any &\hspace{1.4cm}  - &  \hspace{2.3cm} $\pfrac{b_2}{z_4}(z_0) = 0$ \\ 
    %, \pfrac{b_2}{z_4} \not\equiv0 
      &   $a_1 = z_5 $, $a_2(\ol  z_4)$ any & \hspace{1.4cm} -& \\ 
\hline
\end{tabular}
\end{center}
\end{table*}  
    
\section{Discussion of the normal forms}  \label{sec: discussion}
  \vspace{-0.3cm}
% {\it 3.1 \quad Discussion of the normal forms} 
All normal forms are valid around $z_0\in \R^n$, which may be zero or not. Thus all of them can be used around any point (equilibrium or not). 
All forms and the minimal $x$-flat outputs are compatible, that is, {for} a given  $x$-flat system~$\S$ in dimension five,  we can always simultaneously normalize~$\S$ and {a priori given minimal flat output} $\vp$:%, as asserted by:

\begin{proposition}\label{prop compatibility} 
{\it Let~$\S$ be flat at~$x_0$ (at $(x_0, u_0)$ or at $(x_0,$ $u_0,\dot u_0)$) and {$\vp \hspace{-0.05cm}= \hspace{-0.05cm}(\vp_1, \vp_2)$} a minimal $x$-flat output of differential weight at most $n\hspace{-0.05cm}+\hspace{-0.05cm}m\hspace{-0.05cm}+\hspace{-0.05cm}3 \hspace{-0.05cm}= \hspace{-0.05cm}10$ of~$\S$. Then~$\S$~is locally around~$x_0$ static feedback equivalent to one of the forms {$NF^*_1$-$NF^*_{13}$, where $NF^*$ stands for $NF$, $NF'$ or $NF''$,  with $\varphi \hspace{-0.05cm}= \hspace{-0.05cm}(z_1 ,z_2)$ for $NF_1$-$NF_5$, $NF''_8$-$NF''_{12}$ and $NF_{13}$,  $\varphi \hspace{-0.05cm}= \hspace{-0.05cm}(z_1 ,z_3)$ for $NF'_1$-$NF'_6$, and $\varphi \hspace{-0.05cm}= \hspace{-0.05cm}(z_1 ,z_4)$ for $NF_7$, resp.}}
\end{proposition}
  \vspace{-0.2cm}

%\subsection{Discussion of the normal forms $NF^*_1$-$NF^*_7$}  
{\it 4.1\; Discussion of the normal forms $NF^*_1$-$NF^*_7$\quad} 
Theorem~\ref{thm: normal forms dw at most 9} provides normal forms 
{that} are either static feedback linearizable, equivalently, of differential weight $n+m = 7$, {i.e., $\mddiff(\S)= 0$} ({forms} $NF^*_1$), or dynamically linearizable   via a one-fold prolongation of a suitably chosen control, equivalently, of differential weight $n+m +1= 8$, {i.e., $\mddiff(\S)= 1$}   ({forms}  $NF^*_2$,  $NF'_6$ and $NF_7$), or  dynamically linearizable   via a two-fold prolongation of a suitably chosen control, equivalently, of differential weight $n+m +2= 9$, {i.e., $\mddiff(\S)= 2$} ({forms} $NF^*_3$, $NF^*_4$ and $NF^*_5$).  
In the particular case of a static feedback linearizable system, we obtain, as expected the two Brunovsk\'y canonical forms in dimension five consisting of two independent chains of integrators of legths (2,3) for $NF'_1$ and (1,4) for $NF_1$, resp.  
{For the} normal forms describing linearization via a one-fold prolongation 
%
%are particular cases of those of \cite{nicolau2019normal} (that {paper} presents normal forms for flat systems of differential
% weight $n + m + 1$, where~$n$  and~$m$ are arbitrary). In accordance with the results of \cite{nicolau2019normal}, 
%%%%%%%%%%%
one (and only one) nonlinearity is present. 
The normal forms   $NF^*_3$, $NF^*_4$ and $NF^*_5$ corresponding to the differential weight $n+m +2= 9$ 
%(or equivalently, to linearization via a two-fold prolongation)
%% REMOVED 1st submission 
% are particular cases of those presented in \cite{nicolau2020normal}, where  forms for two-input flat systems of differential 
% weight $n + m + 2 = n+4$, {with}~$n$ arbitrary, are given. They 
%%%%%%%%
exhibit (at most) two-nonlinearities. 
Observe that  $NF^*_4$   may {actually} present only one nonlinearity: this happens only if the function~$a_2$ involved in the expression of $\dot z_4$ is identically zero, i.e., the two last equations are given by $\dot z_3 = b_1(\ol z_4) + z_4 v_1$ (the nonlinear one) and $\dot z_4 = z_5$ for $NF_4$ (resp., by $\dot z_3 = b_1(\ol z_4) + (z_4-z_{40}) v_1$  and $\dot z_4 = z_5$ for $NF'_4$).  On the other hand, the normal forms $NF^*_3$ always involve two nonlinear equations: $\dot z_3 = b_1(\ol z_4) + z_4 v_1$ for $NF_3$  or $\dot z_3 = b_1(\ol z_4) + (z_4-z_{40}) v_1$ for $NF'_3$, and $\dot z_4 =  a_1(z) + z_5 v_1$, where the functions~$a_1$ and~$b_1$ may be any, and in particular, both of them can be identically zero; in {the latter case}   $NF^*_3$ are "almost" in the chained form (actually, they are a cascade of a linear subsystem  and the chained form in dimension four). 
The number of nonlinearities in $NF^*_5$ can again vary. 
{Note that all forms $NF^*_1$-$NF^*_7$}  of  Theorem~\ref{thm: normal forms dw at most 9} have a special triangular structure. 
%
% consisting either of two chains of integrators (for the static feedback linearizable case) or of a linear chain and a nonlinear one with at most two nonlinearities (at most two arbitrary functions). We would like to emphasize the fact that {those} arbitrary functions are always consistent with the triangular structure of the system (that is, except {for} the functions $a_i$, they cannot depend on all variables).
% {In fact, all systems exhibit the following structure. 
%
They are cascades of a $k$-dimensional linear system $y_1^{(\rho_1)} = v_1$,  $y_2^{(\rho_2)} = w_2$, where $\rho_1 + \rho_2 = k{\leq 3}$, 
and of a triangular system\vspace{-0.1cm}
\begin{equation}\label{eq: triang cascade}
 \begin{array}{l}
\dot z_j =  d_j(y, \ol z_{j+1}, v_1), \quad \pfrac{d_j}{ z_{j+1}}(z_0)\neq0, \, \dleq {k+1}j 4,\\
\dot z_5 =  v_2\vspace{-0.1cm}
\end{array}
\end{equation}
linked via $w_2 = z_{k+1}$, {where {$d_3 = b$ and $d_4 = a$ of} Theorem~\ref{thm: normal forms dw at most 9}(ii); with $d_3$ absent for $NF^*_2$ and $NF'_6$.}
% {Both} functions $d_j$ are affine with respect to $v_1$ and denoted $d_3 = b_1 + b_2v_1$ (which is 
% absent for $NF^*_2$ and $NF'_6$) and $d_4 = a_1 + a_2v_1$.  
The linear {$y$-}subsystem{s dimension always equals} 2 or 3.

For all normal forms that are not static feedback linearizable,  we see immediately the control to be prolonged. Indeed,  {all of them}
become locally static feedback linearizable after a one- or two-fold prolongation of $v_1$.
Flatness described by Theorem~\ref{thm: normal forms dw at most 9} is local around~$x_0$ {but} for the cases of differential weight $n+m = 7$ or $n+m+1 =8$ (for which the first noninvolutive distribution is~$\mathcal{D}^1$), corresponding to the normal forms $NF^*_1$, $NF'_6$, flatness is global with respect to the control~$u$. On the other hand, in the cases of differential weight $n+m+1=8$ for which $\mathcal{D}^0$ is noninvolutive, and  (for almost all cases) of differential weight $n+m+2 =9$,  the precompensator  always creates singularities in the control space (depending on the state) and flatness  {has} to be considered {around} $(x_0, u_0)$. 
Normal forms $NF^*_5$ are special cases {that exhibit a pretending singularity, where by ``pretending'' we mean that it looks like 
%or pretends to be 
a singularity but it is not; indeed}   
from $\pfrac{(b_1+ b_2 v_{10})}{z_4}(z_0)\neq 0 $,
it follows that
we should be able to normalize either~$b_1$ or~$b_2$ to $z_4$, but $\pfrac{b_2}{z_4}(z_0) = 0$, so we necessarily have $\pfrac{b_1}{z_4}(z_0) \neq 0 $ and it is the function~$b_1$ that has to be normalized to $z_4$ (in the new coordinates, this gives $\pfrac{(b_1+ b_2 v_{10})}{z_4}(z_0) =1+ \pfrac{b_2}{z_4}(z_0) v_{10} = 1 \neq 0 $).
%%%%%%%%% CORRECTION A8
% {Either $\pfrac{b_2}{z_4} \equiv0$ and then the system is flat at any $(z_0, v_0)$ or $\pfrac{b_2}{z_4} \neq0$ and then for any $z_0$ the system is flat in a neighborhood of $(z_0, v_0)$ but that neighborhood can never be taken as $Z_0\times \mathbb{R}^2$, for some open $Z_0\subset Z$.} 
%%%%%%%%%%%%%%

All forms are clearly flat with  $\varphi = (z_1 ,z_2)$ being a
minimal flat output  for $NF_1$-$NF_5$, $\varphi = (z_1 ,z_3)$ for $NF'_1$-$NF'_6$, and $\varphi = (z_1 ,z_4)$ for $NF_7$, resp. 
{An} interesting property of the forms $NF^*_1$-$NF^*_7$ is that each successive time-derivative of the flat output provides the largest  possible number of independent functions of the state~$z$ only (they are actually $x$-maximally flat, see \cite{nicolau2014control} for a formal definition). More precisely, {each successive time-derivative of $\vp_i$, $\dleq 1 i 2$, allows us} to express either two 
(if the control $v_1$ has not been computed yet) or one new state variable 
(if $v_1$ has already been expressed).

%
%% REMOVED 1st submission
% All normal forms of Theorem~\ref{thm: normal forms dw at most 9} can be characterized using geometric conditions either for static feedback linearizability \cite{jakubczyk1980on, hunt1981linear} (forms $NF^*_1$), or dynamic linearizability via a one-fold prolongation \cite{nicolau2016two}
% %, nicolau2017flatness}, 
% % {(with~$\mathcal{D}^1$ the first noninvolutive distribution for $NF'_6$, and $\mathcal{D}^0$ the first noninvolutive distribution for $NF^*_2$ and $NF_7$)}, 
% or   dynamic linearizability via a two-fold prolongation \cite{nicolau2016flatness,gstottner2021necessary}. The normal forms $NF_3$ and $NF_4$ are also particular cases of the triangular forms  characterized in \cite{gstottner2022structurally, gstottner2021structurally}, so the characterizations therein also apply. %\vspace{-0.2cm}
%%%%%%%%%%%

%\subsection{Discussion of the normal forms $NF''_8$-$NF''_{12}$} 
{\it 4.2\; Discussion of the normal forms $NF''_8$-$NF''_{12}$\quad} 
Most {observations made so far} {(triangular structure as cascade of a linear subsystem and a triangular one, $x$-maximally flatness, singularities in the control space involving the control but not its derivatives)} also apply or can be adapted for the
%normal 
forms $NF''_8$-$NF''_{12}$ of  Theorem~\ref{thm: normal forms dw 10}. 
The form $NF''_{8}$ is in the triangular form compatible with the chained form, see \cite{li2016multi, silveira2015flat} for a geometric characterization. 
%Like for $NF_6$, all arbitrary functions of  $NF_{9}$ could be identically zero and then,  $NF_{9}$ coincide with the chained form in dimension five. 
%Most of the forms  $NF''_8$-$NF''_{12}$ fall into a more general class of triangular systems for which we propose a geometric characterization in {a future paper}. %\cite{gstottner2021geometric}. 
{Form 
$NF''_{12}$ coincides with that of~\cite{silveira2007triangular}.} %\vspace{-0.2cm}

%\subsection{Discussion of the normal form $NF_{13}$}\vspace{-0.2cm}
{\it 4.3\; Discussion of the normal form $NF_{13}$\quad}
Similarly to $NF''_8$-$NF''_{12}$, the normal form $NF_{13}$ is also $x$-flat of differential weight $n+m +3= 10$, {i.e., $\mddiff(NF_{13}) = 3$}, with $\vp=(z_1, z_2)$ a minimal flat output, and becomes   dynamically linearizable via a three-fold prolongation of $v_1$, but structurally it is completely different.
First, it does not have a triangular structure and moreover, it is not static feedback equivalent to another form with a triangular structure. 
Second, when expressing the state and control variables in terms of the flat output and its successive derivatives, 
{we have $\mspan\{\md \vp, \md \dot \vp\} \cap \mspan\{\md z\} =\mspan\{\md z_1, \md z_2\}$, where $\mspan\{\md z\} = \mspan\{\md z_1, \ldots, \md z_5\}$. So we are able to express only two states $z_1$ and $z_2$ using four functions $\varphi_1, \varphi_2,\dot \vp_1, \dot \vp_2$ (and we additionally need $\ddot \vp$ to express $z_3$ and $z_4$) implying that the system is not $x$-maximally flat, see \cite{nicolau2014control}.}
Third, the nature of the singularities in the control space, described by relation~\eqref{eq: reg cond}, is also very different since it involves the derivatives of the input $ v$. Condition~\eqref{eq: reg cond} is a consequence of the fact that with $\dot\vp_2$ and $\ddot \vp_2$  we have to express both $z_3$ and $z_4$, so $\mrk \left(\pfrac{(\dot\vp_2, \ddot\vp_2)}{(z_3,z_4)}\right)(z_0, v_0, \dot v_0) = 2$.

Finally, notice that $NF_{13}$ recalls very much  the normal form {of \cite{pomet1997dynamic}} characterizing $(x,u)$-flatness of control systems in dimension four with two inputs: \vspace{-0.1cm}
\begin{equation}\label{eq: normal form pomet}
\begin{array}{lcl lcl}
 \dot z_1 &=& v_1 &\dot z_2 &=& z_4 + z_3v_1\\
 && &\dot z_3 &=& p_0(\ol z_3) + z_4 p_1(\ol z_3) + (q_0(\ol z_3) + z_4 q_1(\ol z_3) ) v_1\\
 &&& \dot z_4 &=& v_2. \vspace{-0.3cm}
\end{array}
\end{equation} 
Indeed, observe first that neither of them is triangular. Second, for both of them, in the coordinates in which $\dot z_2 = z_4 + z_3v_1$ (permute $z_3$ and $z_4$ for $NF_{13}$), the nonlinear functions are affine with respect to the last variable of the nonlinear chain (that is, with respect to $z_4$ for~\eqref{eq: normal form pomet}, resp., to~$z_5$ for $NF_{13}$). Third, when expressing the state {as a} function of $\vp_i^{(j)}$,  the first {time-derivative $\dot \vp$} of $\vp = (\varphi_1, \varphi_2)$ gives no function depending on the state~$z$ {only}. Finally, for both forms, the singularities in the control space depend on the derivatives of the input. \vspace{-0.2cm}

\section{Examples}\label{sec: examples} \vspace{-0.4cm}
{\it Example 5.1.} Consider the following model of an induction motor:\vspace{-0.1cm}
$$
 \begin{array}{lcl c lcl}
  \dot I_d& = &U_d && \dot \omega& = &\mu\psi_d I_q-\tfrac{\tau_L}{J}\\
  \dot I_q& = &U_q&& \dot \psi_d& = &-\eta\psi_d+\eta M I_d\\
  & & & &\dot \rho& = &n_p\omega+\tfrac{\eta M I_q}{\psi_d},
 \end{array}\vspace{-0.1cm}
$$
where {$U_d=u_1$, $U_q=u_2$} are the inputs, see \cite{chiasson1998new} for a detailed explanation of the model. A possible {minimal} $x$-flat output is given by $\vp=(\omega,\rho)$, see e.g. \cite{nicolau2016two} for a derivation of this flat output. This flat output is compatible with $NF'_6$, {in the sense of Proposition~\ref{prop compatibility}}. Indeed, introduce new coordinates $z_1=\vp_1=\omega$, $z_2={L_f\vp_1}=\mu\psi_d I_q-\tfrac{\tau_L}{J}$, $z_3=\vp_2=\rho$, $z_4={L_f\vp_2}=n_p\omega+\tfrac{\eta M I_q}{\psi_d}$ and introduce $v_1={L_f^2\vp_1 + u_2 L_{g_2}L_f\vp_1}=(M I_d- \psi_d)\mu \eta I_q + \mu\psi_d U_q$ by means of an invertible static feedback. We then have
$$
 \begin{array}{lcl c lcl}
  \dot z_1& = &z_2 & &\dot z_3& = &z_4\\
  \dot z_2& = &v_1& &\dot z_4& = &a_1(z_1,z_2,z_4,I_d)+\tfrac{J(z_4-n_pz_1)}{Jz_2+\tau_L}v_1\\
  & & & &\dot I_d& = &U_d.
 \end{array}\vspace{-0.1cm}
$$
Introducing $z_5=a_1(z_1,z_2,z_4,I_d)$ and applying a suitable static feedback then yields $NF'_6$ with $a_2(\ol z_4)=\tfrac{J(z_4-n_pz_1)}{Jz_2+\tau_L}$.\vspace{-0.1cm}
%$$
% \begin{array}{lcl c lcl}
%  \dot z_1& = &z_2 & &\dot z_3& = &z_4\\
%  \dot z_2& = &v_1& &\dot z_4& = &z_5+\tfrac{J(z_4-n_pz_1)}{Jz_2+\tau_L}v_1\\
%  & & & &\dot z_5& = &v_2,
% \end{array}\vspace{-0.1cm}
%$$
%which corresponds to $NF'_6$.\vspace{-0.2cm}

{\it Example 5.2.} As a further example, consider the system
$$
 \begin{array}{lcl c lcl}
  \dot x& = &\cos(\theta_{ 2}-\theta_{ 1})\cos(\theta_{ 1})u_1 && \dot \theta_{ 1}& = &\sin(\theta_{ 2}-\theta_{ 1})u_1\\
  \dot y& = &\cos(\theta_{ 2}-\theta_{ 1})\sin(\theta_{ 1})u_1&& \dot \theta_{ 2}& = & \omega_2\\
  & & & &\dot {\omega}_2& = & u_2,
 \end{array}\vspace{-0.1cm}
$$
which, {for a nonholonomic car}, models the kinematics, {see, e.g., \cite{fliess61vine},} and the steering dynamics {for the angular acceleration $\ddot \theta_2 = u_2$}. A minimal $x$-flat output of this system {can be taken as} $\vp=(x,y)$. This flat output is compatible, {in the sense of Proposition~\ref{prop compatibility}}, with $NF_9''$. For transforming the system into $NF_9''$, introduce new coordinates $z_1=\vp_1=x$, $z_2=\vp_2=y$ and {set} $v_1={u_1L_{g_1}}\vp_1=\cos(\theta_{ 2}-\theta_{ 1})\cos(\theta_{ 1})u_1$, resulting in
$$
\begin{array}{rcl rcl}
  \dot z_1& = &v_1 & \dot z_2& = &\tan(\theta_{ 1})v_1 , \\
  & & & \dot \theta_{ 1}& = &\tfrac{\tan(\theta_{ 2}-\theta_{ 1})}{\cos(\theta_{ 1})} v_1,\\
  & & & \dot \theta_{ 2}& = & \omega_2,\\
  & & & \dot {\omega}_2& = & u_2.
 \end{array}
$$
Next, successively introduce the functions multiplying $v_1$ as new states, i.e., {set} $z_3=\tan(\theta_{ 1})$ and after that, {define}~$z_4$ in a similar fashion, resulting in $\dot z_2=z_3v_1$, $\dot z_3=z_4v_1$, ${\dot z_4}=a_1(z_3,z_4, \omega_2)+a_2(z_3,z_4)v_1$. Finally, introducing $z_5=a_1(z_3,z_4, \omega_2)$ and applying a suitable static feedback yield a representation in $NF_9''$.\vspace{-0.2cm}

\begin{small}
\bibliographystyle{plain}
\bibliography{bibliographie}
\end{small}

\end{document}